\documentclass[11pt]{article}
\usepackage[a4paper,top=3cm,bottom=2cm,left=3cm,right=3cm,marginparwidth=1.75cm]{geometry}
\usepackage{graphicx}%
\usepackage{multirow}%
\usepackage{amsmath,amssymb,latexsym,amstext,amsthm}
\usepackage{amsfonts}%
\usepackage{mathrsfs}%
\usepackage{hyperref}

\theoremstyle{thmstyleone}%
\newtheorem{theorem}{Theorem}[section]        
\newtheorem{lemma}{Lemma}[section]
\newtheorem{corollary}{Corollary}[section]
\newtheorem{definition}{Definition}[section]

\theoremstyle{thmstyletwo}%
\newtheorem{remark}{Remark}[section]%

\newcommand{\Pochh}[2]{\left( #1 \right)_{#2}}
\newcommand{\supp}{\text{supp }}
\newcommand{\spn}{\text{span }}

\providecommand{\keywords}[1]
{
  \small	
  \textbf{\textit{Keywords:}} #1
}

\providecommand{\subjclass}[1]
{\small	
  \textbf{\textit{2010 Mathematical Subject Classification:}} #1
}

\title{Nearest neighbor recurrence relations for Meixner-Angelesco multiple orthogonal polynomials of the second kind}

\author{Jorge Arves\'u Carballo\\
\small Department of Mathematics, Universidad Carlos III de Madrid, Avda. de la
Universidad, 30 \\
\small Legan\'es, 28911 Madrid, Spain\\
Alejandro J. Quintero Roba\\
\small Department of Mathematics, Baylor University, 1410 S
4th Street. Waco, 76706, Texas, USA}

\date{October 1, 2023}

\begin{document}
\maketitle

\begin{abstract} This paper studies a new family of Angelesco multiple orthogonal polynomials with shared orthogonality conditions with respect to a system of weight functions, which are complex analogues of Pascal distributions on a legged star-like set. The emphasis is placed on the algebraic properties such as the raising operators, the Rodrigues-type formula, the explicit expression of the polynomials, and the nearest neighbor recurrence relations.
\end{abstract}

\vspace{.2cm}

\noindent\subjclass{33C45, 33C47, 41A28, 42C05}\\
\noindent\keywords{Angelesco polynomials; Classical orthogonal polynomials; Discrete orthogonal polynomials; Multiple orthogonal polynomials; Meixner polynomials; Ro\-dri\-gues-type formula; Recurrence relations}

\section{Introduction}\label{intro}

This article focuses on the algebraic properties of a new family of type II multiple orthogonal polynomials, which extend the notion of discrete orthogonal polynomials on the real line \cite{bib:nikiforov-uvarov-suslov}. Multiple orthogonal polynomials are a natural extension of orthogonal polynomials. They are linked to the simultaneous rational approximation of a system of analytic functions \cite{bib:nikishin-sorokin,bib:Sorokin3}. Hermite's proof of the transcendence of the number $e$ \cite{Hermite_exp} use this notion. Depending on the simultaneous approximation problem, some types of multiple orthogonal polynomials appear (type I, type II, and mixed type). See \cite[Chapter 4]{bib:nikishin-sorokin} as well as \cite{bib:aptekarev,bib:Brezinski-Iseghem,bib:Bruin,bib:Bruin2} for simultaneous rational approximation (also known as Hermite-Pad\'e approximation) and types of multiple orthogonal polynomials.

In \cite{bib:arvesu-coussement-vanassche} some type II discrete multiple orthogonal polynomials were investigated. In particular, the multiple Meixner polynomials which verify orthogonality relations with respect to $r>1$ different Pascal distributions (negative binomial distributions). Two kinds of multiple Meixner polynomials that are distinguished by the selection of parameters in the Pascal distribution were addressed. In \cite{bib:miki-tsujimoto-vinet-zhedanov1} was given a mathematical approach to a physical phenomenon involving a chain of Hamiltonians and the multiple Meixner polynomials of the first kind. Later in \cite{bib:ndayiragije1-vanassche}, the authors constructed $r$ non-Hermitian oscillator Hamiltonians that are simultaneously diagonalizable and for which the common eigenstates are expressed in terms of multiple Meixner polynomials of the second kind. In this paper we study a different type of  multiple Meixner polynomials, which are orthogonal with respect to $r$ complex analogues of the Pascal distributions on an $r$-star and study their algebraic properties. Among these properties, special attention is given to the nearest neighbor recurrence relations. This question has been addressed in a few recent papers dealing with multiple orthogonal polynomials on an $r$-star \cite{bib:arvesu-quintero1,bib:leurs-vanassche2,bib:leurs-vanassche1}. Knowing the coefficients of the recurrence relations one can study the weak asymptotics (see \cite{bib:aptekarev-arvesu} for the multiple Meixner polynomials of the first and second kind on the real line). Also the ratio asymptotic behavior of multiple orthogonal polynomials. Moreover, having such recurrence relations one can construct Christoffel-Darboux kernels, which plays important role in correlation kernel, for instance in the unitary random matrix model with external source (see \cite{bib:kuijlaars-congress} for the connection of multiple orthogonal polynomials and random matrix theory).

The structure of the paper is as follows. In Section \ref{sec:section2} we provide the background material. Sections \ref{sec:meix2_examples} and \ref{recu} contain the main results. First, we define the system of measures used in the paper and study a new family of type II discrete Angelesco multiple orthogonal polynomials on an $r$-star, namely, the Meixner-Angelesco polynomials of the second kind. Then, in subsection \ref{raising} we obtain the raising operators and state their commutative property. Base on these operators the Rodrigues-type formula is obtained. Section \ref{recu} deals both with the explicit expression of the polynomials and the deduction of the nearest neighbor recurrence relation. The implemented procedure in the study of the recurrence relations is similar to the one used in \cite{bib:arvesu-ramirez1}, which differs from those used in \cite{bib:arvesu-coussement-vanassche,bib:arvesu-quintero1,bib:vanassche}. The explicit series representation of the polynomials is given. The paper ends with some concluding remarks in Section \ref{conclusions}.

\section{Background material}\label{sec:section2}

The orthogonal polynomials on the real line $P_n(x)$, $n\in\mathbb{Z}_{+}=\{0,1,2,\ldots\}$, with respect to a positive measure $\mu$ (with finite moments) are such that each polynomial has $\deg P_n=n$ and 
satisfy the following orthogonality relations 
\begin{equation}\label{ortho_classic_eq_1}
\displaystyle\int_{\Omega\subset\mathbb{R}}P_n(x)x^k\, \text{d}\mu (x)=0,\qquad k=0,\dots,n-1.
\end{equation}
From \eqref{ortho_classic_eq_1} the polynomial $P_n$ is determined up to a multiplicative factor.
For discrete orthogonal polynomials \cite{bib:nikiforov-uvarov-suslov}, the measure $\mu$ is discrete 
\begin{equation*}
\mu=\sum_{k=0}^{N}\rho_k\delta_{x_k},\quad \rho_k>0,\quad x_k\in\mathbb{R},\quad \mbox{and}\ N\in\mathbb{Z}_{+}\cup \{+\infty\},
\end{equation*}
where $\delta_{x_k}$ denotes the Dirac measures on the $N+1$ points $x_0,\ldots,x_N$ and $\rho_k$ are the weights.  These orthogonality relations on the linear lattice $x_k = k$, $k=0,1,\ldots,N$ can be expressed as
\begin{equation}\label{ortho_classic_eq_2}\sum \limits_{k=0}^{N}p_n(k)(-k)_j\, \rho_k=0, \qquad j=0,\ldots,n-1,
\end{equation}
where $(x)_j=x(x+1)\dots(x+j-1)$, if $j>0$ and $(x)_{0}=1$, is the Pochhammer symbol. Next we take $\rho_k=\rho(k)$ (see Table \ref{table1}).

The classical discrete orthogonal polynomials on a linear lattice are those of Charlier, Meixner, Kravchuk, and Hahn.  The corresponding weight functions $\rho(x)$ for the discrete measures $\mu$ with unbounded support, i.e. $N=\infty$, are given in the following table (see \cite{bib:nikiforov-uvarov-suslov}).
\begin{table}[!h]
\begin{tabular}
{|c|c|c|}
\hline
Weight & Charlier $\rho^{a}(x)$ & Meixner $\rho^{\beta,c}(x)$ \\
\hline\hline
$\rho(x)$ & $\dfrac{a^x}{\Gamma(x+1)}$ & $\dfrac{c^x\Gamma(\beta +x)}{\Gamma(x+1)}$ \\
\hline
\end{tabular}
\caption{Weight functions for Charlier and Meixner polynomials.}
\label{table1}
\end{table}

For the Charlier polynomials $C_n^{a}(x)$ 
the weight function $\rho^{a}(x)$, $x\in \mathbb{Z}_{+}$, $a>0$, is a Poisson distribution on $\mathbb{Z}_{+}$, and the orthogonality relations \eqref{ortho_classic_eq_2} are:
\[\sum \limits_{k=0}^{+\infty}C_n^{a}(k)(-k)_j\frac{a^k}{k!}=0, \qquad  j=0,\ldots,n-1.\]

For the Meixner polynomials $M_n^{\beta,c}(x)$ (with $\beta>0$ and $0<c<1$) the weight function $\rho^{\beta,c}(x)$,
$x\in \mathbb{Z}_{+}$, is a negative binomial distribution (Pascal distribution) on $\mathbb{Z}_{+}$, and the orthogonality relations \eqref{ortho_classic_eq_2} are:
\[\sum \limits_{k=0}^{+\infty}M_n^{\beta,c}(k)(-k)_j\frac{(\beta)_k}{k!}c^k
=0, \qquad j=0,\ldots,n-1.\]

The weight functions for classical discrete orthogonal polynomials satisfy a first order difference equation (Pearson's equation) $\bigtriangleup \Bigl(\sigma(x)\rho(x)\Bigr)= \tau(x)\rho(x)$, where $\sigma$ is a polynomial of degree at most $2$, $\tau$ a polynomial of degree $1$, and 
$\bigtriangleup f(x)=f(x+1)-f(x)$, the forward difference of $f$ in $x$. The backward difference of $f$ in $x$ is 
\begin{equation}\label{triang_nabla}
\bigtriangledown f(x)=f(x)-f(x-1).
\end{equation}
For the above $P_{n}(x)$, the polynomials $\sigma(x)$ and $\tau(x)$ are given in the following table \cite{bib:arvesu-coussement-vanassche}, \cite{bib:nikiforov-uvarov-suslov}:
\begin{table}[!h]
\begin{tabular}
{|c|c|c|}
\hline
$P_{n}(x)$ & $C_n^{a}(x)$ & $M_n^{\beta,c}(x)$ \\
\hline\hline
$\sigma(x)$ & $x$ & $x$ \\
\hline
$\tau(x)$ & $a-x$ & $(c-1)x+\beta c$ \\
\hline
\end{tabular}
\caption{Polynomials $\sigma(x)$ and $\tau(x)$ in the Pearson's equation for Charlier and Meixner polynomials.}
\label{tab1-1}
\end{table}

Each of the classical (monic) discrete orthogonal polynomials has a
raising operator, which can be derived from the orthogonality relation \eqref{ortho_classic_eq_2} and the use of summation
by parts. These operators are
\begin{align}
&-\left(\dfrac{\Gamma(x+1)}{a^{x-1}}\bigtriangledown  \frac{a^x}{\Gamma(x+1)}\right)C_n^{a}(x) 
 = C_{n+1}^{a}(x),\notag\\
&-\left(\dfrac{\Gamma(x+1)}{(\beta)_{x-1}c^{x-1}(1-c)}
\bigtriangledown  \frac{(\beta)_x}{\Gamma(x+1)}c^x\right)M_n^{\beta,c}(x) 
 =  M_{n+1}^{\beta-1,c}(x).\label{meixner_classic}
\end{align}

Using the raising operator several times, one gets a
\textit{Rodrigues-type formula} for the polynomials \cite{bib:arvesu-coussement-vanassche}. From this formula, one can obtains an explicit expression for the classical (monic) discrete orthogonal polynomials in terms of hypergeometric functions \cite{bib:koekoek-lesky-swarttouw}
\begin{equation}\begin{array}{rcl}
C_n^{a}(x)& = & (-a)^n \ _2F_0
\left(
\left.
\begin{array}{c}
-n,-x \\
-
\end{array}
\right|
-\frac{1}{a} \right),\\
M_n^{\beta,c}(x)& = & \dfrac{c^{n}}{(c-1)^n}(\beta)_n \ _2F_1
\left(
\left.
\begin{array}{c}
-n,-x \\
\beta
\end{array}
\right|
1-\frac{1}{c} \right).
\end{array}
\label{hyper_classic}
\end{equation}

These classical (monic) discrete orthogonal polynomials $P_n(x)$ satisfy a three-term recurrence relation 
\[xP_n(x)=P_{n+1}(x)+b_nP_n(x)+a_nP_{n-1}(x),\quad n \ge 0, \]
with $a_n > 0$ and initial conditions $P_0=1$  and  $P_{-1}=0$. 

Substituting the explicit expression given in \eqref{hyper_classic} in the recurrence relation, one can compute the above recurrence coefficients by comparing coefficients \cite{bib:arvesu-coussement-vanassche}. Indeed, they are given in table \ref{tab1-2}.
\begin{table}[!h]
\begin{tabular}
{|c|c|c|}
\hline
$P_{n}(x)$ & $a_n$ & $b_{n}$ \\ 
\hline\hline
$C_n^{a}(x)$ & $an$ & $a+n$ \\
\hline
$M_n^{\beta,c}(x)$ & $\frac{cn(\beta+n-1)}{(1-c)^2}$ & $\frac{n+(\beta+n)c}{1-c}$ \\
\hline
\end{tabular}
\caption{Recurrence coefficients for Charlier and Meixner polynomials.}
\label{tab1-2}
\end{table}

An extension of classical orthogonal polynomials are the multiple orthogonal polynomials (also called Hermite-Pad\'e polynomials). They satisfy orthogonality relations shared with respect to a system of measures and are linked with the simultaneous rational approximation of a system of analytic functions~\cite{bib:nikishin-sorokin}. Below, we address two of these systems of measures: The AT-system and Anglesco system. In \cite{bib:arvesu-coussement-vanassche} some Hermite-Pad\'e polynomials orthogonal with respect to an AT-system of discrete measures supported on the linear lattice were studied.
An AT system ~\cite{bib:arvesu-coussement-vanassche,bib:nikishin-sorokin} of $r$ positive discrete measures consists of a set of measures
\[\mu_j=\sum_{k=0}^{N} \rho_{j,k} \delta_{x_{k}},\qquad \rho_{j,k}>0,\
x_k\in\mathbb{R},\ \ N\in
\mathbb{N} \cup \{+\infty\},\ \ \ j=0,\ldots,r-1,
\]
on $\Omega_{j}\subset\mathbb{R}$, so that $\mbox{supp}(\mu_j)$ is the closure of $\{x_k\}_{k=0}^N$ and
$\Omega_j=\Omega$ for each $j=0,\ldots,r-1$.
Moreover, $\rho_{j,k}=\rho_j(x_k)$, $k=1,\ldots,N$, $j=0,\ldots,r-1$, for weight functions $\rho_0(x),\ldots,\rho_{r-1}(x)$ such that the following system  
\[\begin{array}{c}
\rho_0(x),x\rho_0(x),\ldots,x^{n_0-1}\rho_0(x),\ldots,
\rho_{r-1}(x),x\rho_{r-1}(x),\ldots,x^{n_{r-1}-1}\rho_{r-1}(x)
\end{array}
\]
be a Chebyshev system on $\Omega$ for each multi-index $\vec{n}=(n_0,\ldots,n_{r-1})\in\mathbb{Z}_{+}^{r}$ with 
length $|\vec{n}|=n_0+\cdots+n_{r-1} < N+1$. For this system of measures on the linear lattice, define 
the type II discrete multiple orthogonal polynomial of $\deg P_{\vec{n}}=|\vec{n}|$ that satisfies the orthogonality conditions
\begin{equation}\label{type2_mdop}
\sum \limits_{k=0}^{N}P_{\vec{n}}(k)(-k)_j\,\rho_{\ell}(k)=0,\quad j=0,\ldots,n_{\ell}-1,\quad
\ell=0,\ldots,r-1.
\end{equation}

The discrete multiple orthogonal polynomials studied in \cite{bib:arvesu-coussement-vanassche} are those of Multiple Charlier, Multiple Meixner of the first and second kind, Multiple Kravchuk, and Multiple Hahn. The corresponding components of the vector weight functions for measures with unbounded support ($N=\infty$) are given in the following table.
\begin{table}[!h]
\begin{tabular}
{|c|c|c|}
\hline
Charlier $\rho^{a_{\ell}}(x)$ & Meixner I $\rho^{\beta,c_{\ell}}(x)$ & Meixner II $\rho^{\beta_{\ell},c}(x)$ \\
\hline\hline
$\dfrac{a_{\ell}^x}{\Gamma(x+1)}$ & $\dfrac{c_{\ell}^x\Gamma(\beta +x)}{\Gamma(x+1)}$ & 
$\dfrac{c^x\Gamma(\beta_{\ell} +x)}{\Gamma(x+1)}$
\\
\hline
\end{tabular}
\caption{Weight functions $\rho_{\ell}(x)$ involved in the orthogonality conditions \eqref{type2_mdop}.}
\label{tab1}
\end{table}

For the multiple Charlier polynomials $C_{\vec{n}}^{\vec{a}}(x)$ with the set of different positive parameters $(a_0,\ldots,a_{r-1})=\vec{a}$, the orthogonality relations \eqref{type2_mdop} are:
\begin{equation*}
\sum_{k=0}^{+\infty} C_{\vec{n}}^{\vec{a}}(k)(-k)_j\dfrac{a_{\ell}^k}{\Gamma(k+1)}= 0,\quad j=0,\ldots,n_{\ell}-1,\quad
\ell=0,\ldots,r-1.
\end{equation*}

For the multiple Meixner polynomials of the first kind $M_{\vec{n}}^{\beta,\vec{c}}(x)$
with different parameters $(c_0,\ldots,c_{r-1})=\vec{c}$, $0<c_{\ell}<1$, and $\beta \neq 0,-1,-2,\ldots$, the orthogonality relations \eqref{type2_mdop} are:
\begin{equation*}
\sum_{k=0}^{+\infty}M_{\vec{n}}^{\beta,\vec{c}}(k)(-k)_j\dfrac{c_{\ell}^k\Gamma(\beta +k)}{\Gamma(k+1)}=0,
\quad j=0,\ldots,n_{\ell}-1,\quad
\ell=0,\ldots,r-1.
\end{equation*}

For the multiple Meixner polynomials of the second kind $M_{\vec{n}}^{\vec{\beta},c}(x)$
with parameters $(\beta_0,\ldots,\beta_{r-1})=\vec{\beta}$ and $c$, the orthogonality relations \eqref{type2_mdop} are:
\begin{equation*}
\sum_{k=0}^{+\infty}M_{\vec{n}}^{\vec{\beta},c}(k)(-k)_j\dfrac{c^k\Gamma(\beta_{\ell} +k)}{\Gamma(k+1)}=0,
\quad j=0,\ldots,n_{\ell}-1,\quad
\ell=0,\ldots,r-1,
\end{equation*}
with $0<c<1$, $\beta_\ell \neq 0,-1,-2,\ldots$ and pairwise different.

Using summation by parts in the above orthogonality relations one can find the raising operations
\begin{align*}
&-\frac{a_{\ell}}{\rho^{a_{\ell}}(x)}\bigtriangledown \left( \rho^{a_{\ell}}(x)C_{\vec{n}}^{\vec{a}}(x) \right)=
C_{\vec{n}+\vec{e}_{\ell}}^{\vec{a}}(x),\\
&-\frac{c_{\ell}}{1-c_{\ell}}\frac{1}{\rho^{\beta-1,c_{\ell}}(x)}\bigtriangledown \left( \rho^{\beta,c_{\ell}}(x)M_{\vec{n}}^{\beta,\vec{c}}(x)\right)=
M_{\vec{n}+\vec{e}_{\ell}}^{\beta-1,\vec{c}}(x),\\ 
&-\frac{c}{1-c}\frac{1}{\rho^{\beta_{\ell}-1,c}(x)}\bigtriangledown \left( \rho^{\beta_{\ell},c}(x)M_{\vec{n}}^{\vec{\beta},c}(x)\right)=
M_{\vec{n}+\vec{e}_{\ell}}^{\vec{\beta}-\vec{e}_{\ell},c}(x).\\
\end{align*}
where $\ell=0,\ldots,r-1$ and $\vec{e}_{\ell}$ denotes the standard $r$-dimensional unit vector with the $\ell$-th entry equals $1$ and $0$ otherwise.

A repeated application of these raising operators gives the following Rodrigues-type formulas
\begin{align} 
C_{\vec{n}}^{\vec{a}}(x)&=(-1)^{|\vec{n}|}\left(\prod_{\ell=0}^{r-1}a_{\ell}^{n_{\ell}}\right)\Gamma(x+1)\prod_{j=0}^{r-1} \left(\frac{1}{a_{j}^x}\bigtriangledown^{n_{j}}a_{j}^x\right)\frac{1}{\Gamma(x+1)},\notag\\
M_{\vec{n}}^{\beta,\vec{c}}(x)&=\prod_{\ell=0}^{r-1}
\left(\frac{c_{\ell}}{c_{\ell}-1}\right)^{n_{\ell}}
\frac{\Gamma(x+1)}{\Gamma(\beta+x)}\prod_{j=0}^{r-1}\left(\frac{1}{c_j^x}\bigtriangledown^{n_j}c_j^x\right)
\frac{\Gamma(\beta+|\vec{n}|+x)}{\Gamma(x+1)},\notag\\
M_{\vec{n}}^{\vec{\beta},c} (x) &= \left( \dfrac{c}{c-1} \right)^{|\vec{n}|} \dfrac{\Gamma(x+1) }{c^{x}}\prod\limits_{j=0}^{r-1} \left(\Gamma(x+\beta_j)^{-1}\bigtriangledown ^{n_j}  \Gamma(x+ n_j+\beta_j)\right) \dfrac{c^x}{\Gamma(x+1)}.\label{eq:m2_mclassic}
\end{align}
From the Rodrigues-type formulas one can find an explicit expression of these monic polynomials
$$
P_{\vec{n}}(x) = \sum\limits_{m=0}^{\lvert\vec{n}\rvert} \alpha_{m}^{\vec{n}} x^{m}.
$$
Moreover, these polynomials satisfy a nearest neighbor recurrence relations \cite{bib:vanassche}

\begin{equation*}\label{eq:recurrence_relation}
\left(x - b_{\vec{n},k} \right)P_{\vec{n}}(x)  - P_{\vec{n}+\vec{e}_k}(x) = \sum_{\ell=0}^{r-1} d_{\vec{n},\ell}  P_{\vec{n}-\vec{e}_\ell}(x), 
\quad k=0,\ldots, r-1,
\end{equation*}
\ \\
where 
\begin{align*}\label{eq:rr_linear_combination1}
b_{\vec{n},k}= \alpha_{\lvert\vec{n}\lvert-1}^{\vec{n}} - \alpha_{\rvert\vec{n}\lvert}^{\vec{n}+\vec{e}_k}, \quad 
d_{\vec{n},\ell} = \dfrac{\displaystyle\int x P_{\vec{n}}(x)\Pochh{-x}{n_{\ell}-1}\text{d} \mu_\ell}{\displaystyle\int  P_{\vec{n}-\vec{e}_{\ell}}(x)\Pochh{-x}{n_{\ell}-1}\text{d} \mu_\ell}, \quad \ell=0,\ldots,r-1. \notag
\end{align*}

Now, we address the Angelesco system of $r$ measures \cite{bib:angelesco,bib:nikishin-sorokin}. This system consists of orthogonality measures $\mu_{\ell}$ $(\ell=0,\ldots,r-1)$ supported on $\Omega_{\ell}$, where $\Omega_{\ell}^{\circ}=\Omega_{\ell}\setminus\{0\}$ and $\Omega_i^{\circ} \cap \Omega_j^{\circ} = \emptyset$,  for  $i\neq j$ (see also \cite{bib:sorokin2}). In \cite{bib:arvesu-quintero1}, an Angelesco system of discrete measures were considered. The derived polynomials (Angelesco multiple orthogonal polynomials) involve orthogonality relations distributed over $r$ discrete complex measures (with finite moments) supported on an $r$-star defined by the intervals $[0,e^{2\pi i k/r}\infty)$, $k=0,\ldots, r-1$. Two type II discrete multiple orthogonal polynomials on an $r$-star, namely, Charlier-Angelesco and Meixner-Angelesco polynomials of the first kind were investigated. The location of zeros as well as the recurrence relations were addressed.  

In the sequel we assume that $r\in \mathbb{N}$ and $\omega^j$ is the $r$-th roots of $1$, that is $\omega = e^{2\pi i/r}$, $j=0,\ldots, r-1$. Denote by $\Omega_0,\ldots, \Omega_{r-1}$ the rays on the complex plane given as counterclockwise rotations of the positive real axis, i.e.
$$
\Omega_j =  \left[ 0, \omega^j\infty \right) = \omega^j \cdot \mathbb{R_+}\cup \{ 0 \}.
$$ 
These rays generate the $r$-star on the complex plane $\mathbb{C}$
\begin{equation}\label{r-star}
\bigcup\limits_{\ell=0}^{r-1}\Omega_\ell.
\end{equation}

On each $\Omega_j$, $j=0,\ldots,r-1$, one considers a continuous function $\rho_j$ (weight function) and a discrete measure on the mass points $z_{j,0},\ldots,z_{j,N}\in \Omega_{j}$,
\begin{equation} \label{eq:particular_measures}
\mu_j 
=\sum_{k=0}^{N} \rho_{j}\delta_{z_{j,k}},\quad N \in \mathbb{N} \cup \{ \infty \},
\end{equation}
with finite moments and $\supp {\mu_j} = \Omega_j$, with $\Omega_i^{\circ} \cap \Omega_j^{\circ} = \emptyset$,  for  $i\neq j$.

\begin{definition}[\cite{bib:arvesu-quintero1}] \label{def:typeII}
The type II discrete-Angelesco polynomial $P_{\vec{n}}$ for the multi-index $\vec{n}$ on the $r$-star \eqref{r-star} is the polynomial of degree $\leq r\lvert\vec{n}\rvert$ defined by the orthogonality relations
\begin{align}\label{eq:discrete_orth_conditions}
\displaystyle\sum_{k=0}^{\infty} P_{\vec{n}}(z_{\ell,k}) \Pochh{-z_{\ell,k}^r}{j} \rho_{\ell}(z_{\ell,k}) &= 0,\quad 0\leq \ell  \leq r-1,\quad 0\leq j\leq n_{\ell}-1.
\end{align} 
\end{definition}
Here the Pochhammer symbol $\Pochh{-z^r}{n_{\ell}-1}$ is a polynomial of degree $r(n_{\ell}-1)$ and $P_{\vec{n}}$ is a polynomial in powers of $z^r$. 

The orthogonality relations \eqref{eq:discrete_orth_conditions} give a linear algebraic systems of $\lvert\vec{n}\rvert$ equations for the $\lvert\vec{n}\rvert+1$ coefficients of $P_{\vec{n}}$, with matrix 
\begin{eqnarray}
\label{matrix}
\mathcal{M} & = & \left(
\begin{array}{c}
\mathcal{M}_0(n_1)\\
\vdots\\
\mathcal{M}_{r-1}(n_r)
\end{array} \right),\quad \mathcal{M}_j(n_{j+1}) = \left(
\begin{array}{cccc}
m_0^{(j)} & m_1^{(j)} & \ldots & m_{\lvert\vec{n}\rvert}^{(j)}\\
\vdots & \vdots & & \vdots\\
m_{n_j-1}^{(j)} & m_{n_j}^{(j)} & \ldots &
m_{\lvert\vec{n}\rvert+n_j-1}^{(j)}
\end{array} \right),
\end{eqnarray}
where $\mathcal{M}_j(n_{j+1})$ is the $n_{j}\times(\lvert{\vec{n}}\rvert+1)$ matrix of moments 
\begin{equation}\label{eq:particular_moments}
m_{\ell}^{(j)} = \int_{0}^{\omega^{j}\infty} (z^r)^{\ell}\text{d}\mu_j = \sum_{k=0}^{N} \left( z_{j,k}^r \right)^{\ell} \rho_j (z_{j,k}^r),\quad j=0,\ldots,r-1,\quad \ell = 0,1, \ldots
\end{equation}
This linear algebraic system has always a nontrivial solution. We will deal with those systems, whose solution $P_{\vec{n}}$ for all indices $\vec{n}\in\mathbb{N}^{r}$ is unique up to a nonzero multiplicative factor and also of degree exactly $\lvert\vec{n}\rvert$ (the monic polynomial exists and will be unique) i.e., the system of measures \eqref{eq:particular_measures} is a perfect system \cite{bib:Mahler} (all indices are normal). This happen when the matrix 
$\mathcal{N}$ given from $\mathcal{M}$ by deleting the last column has rank $\lvert\vec{n}\rvert$. Hence, 
\begin{equation}
\displaystyle P_{\vec{n}}(z) = \frac{1}{\det \mathcal{N}} \det\left(
\begin{array}{c}
\mathcal{M}
\\ 
\hline
\mathcal{L}
\end{array} \right),\quad \mathcal{L}=\left(\begin{array}{cccc}1 & z^{r} & \cdots &
z^{r\lvert\vec{n}\rvert}\end{array}\right).\label{Pn_matrix_eq}
\end{equation}
For the coefficients of the monic $P_{\vec{n}}$ we use the following notation
\begin{equation*}
P_{\vec{n}}(z) = \sum\limits_{m=0}^{\lvert\vec{n}\rvert} \alpha_{m}^{\vec{n}} z^{rm},
\quad 
 P_{\vec{n}+\vec{e}_k}(z) = \sum\limits_{m=0}^{\lvert\vec{n}\rvert+1} \alpha_{m}^{\vec{n}+\vec{e}_k} z^{rm},\quad \mbox{where}\quad \alpha_{\lvert\vec{n}\rvert}^{\vec{n}} 
= \alpha_{\lvert\vec{n}\rvert+1}^{\vec{n}+\vec{e}_k} = 1.
\end{equation*}
 
Now, we address some properties of these polynomials, namely, the zero location theorem and the recurrence relation theorem (also valid for the new Meixner-Angelesco polynomials of the second kind). For the first property, we need the notion of $\omega$-symmetry \cite{bib:arvesu-quintero1}.

\begin{definition}[\cite{bib:arvesu-quintero1}] \label{def:w-symmetric}
The system of weight functions $\{ \rho_{0}(z),\rho_1(z)\ldots, \rho_{r-1}(z)\}$ is said to be $\omega$-symmetric if
for every $j=0,\ldots,r-1$ and $k\in\mathbb{N}$
\begin{eqnarray*}
\rho_j(\omega^{k}z)=\rho_j(z).
\end{eqnarray*}
\end{definition}
Since the measures \eqref{eq:particular_measures} are defined via the weight functions $\rho_j$ ($j=0,\ldots,r-1$), this $\omega$-symmetry allows to interpret the orthogonality relation \eqref{eq:discrete_orth_conditions} for $j=1,\ldots,r-1$, as rotated copies of the relation with real support 
$$\displaystyle\int_{0}^{\infty} (z^r)^k P_{\vec{n}} (z) \text{d} \mu_0 = 0, \quad k = 0, \ldots, n_0-1.
$$

The following theorems involving the type II discrete Angelesco multiple orthogonal polynomials on the $r$-star were proved in \cite{bib:arvesu-quintero1}. 

\begin{theorem}[\cite{bib:arvesu-quintero1}] \label{location_zeros} Let $\vec{n}=(n,n,\ldots,n)$ and $\rho_{j}(z)$ be the $\omega$-symmetric weight functions supported on $\Omega_j$ $(j=0,\ldots,r-1)$ of the $r$-star \eqref{r-star}. Then, all the zeros of $P_{\vec{n}}$ in \eqref{Pn_matrix_eq} are simple and lie in the rays of the $r$-star.
\end{theorem}

\begin{theorem}[\cite{bib:arvesu-quintero1}] \label{theorem:recurrence_relation}
Suppose all the multi-indices $\vec{n}\in\mathbb{N}^{r}$ are normal for the system of measures $\mu_0, \ldots, \mu_{r-1}$. Then, the corresponding type II discrete Angelesco multiple orthogonal polynomials on the $r$-star satisfy the recurrence relation
\begin{equation*}
\left(z^r - b_{\vec{n},k} \right)P_{\vec{n}}(z)  - P_{\vec{n}+\vec{e}_k}(z) = \sum_{\ell=0}^{r-1} d_{\vec{n},\ell}  P_{\vec{n}-\vec{e}_\ell}(z), 
\quad k=0,\ldots, r-1,
\end{equation*}
\ \\
where 
\begin{align*}\label{eq:rr_linear_combination1}
b_{\vec{n},k}= \alpha_{\lvert\vec{n}\lvert-1}^{\vec{n}} - \alpha_{\rvert\vec{n}\lvert}^{\vec{n}+\vec{e}_k}, 
\quad d_{\vec{n},\ell} = \dfrac{\displaystyle{\int_{\Omega_{\ell}}} z^r P_{\vec{n}}(z)\Pochh{-z^r}{n_{\ell}-1}\text{d} \mu_\ell}{\displaystyle{\int_{\Omega_{\ell}}} P_{\vec{n}-\vec{e}_{\ell}}(z)\Pochh{-z^r}{n_{\ell}-1}\text{d} \mu_\ell}, \quad \ell=0,\ldots,r-1. \notag
\end{align*}
\end{theorem}

The explicit expression of these recurrence coefficients for Charlier-Angelesco and Meixner-Angelesco polynomials of the first kind are given in \cite{bib:arvesu-quintero1}.

In the sequel, we will use the following relations involving the backward difference operator $\nabla$ acting on composite functions 
\begin{equation}
\nabla f(s(z)) = f(s(z)) - f(s(z)-1),\quad 
\nabla f(s(z)+1)=\Delta f(s(z)),\label{relation_delta_nabla}
\end{equation} with $s(z)\pm 1$ also in the domain of $f$.
Moreover, 
\begin{align}
\nabla f(s(z))g(s(z)) &= f(s(z)) \nabla g(s(z)) + g(s(z)-1) \nabla f(s(z)),\notag\\
\nabla^n f(s(z)) &= \displaystyle{\sum\limits_{k=0}^{n} {n\choose{k}}}\left(-1\right)^k f(s(z)-k),\label{nth_nabla}\\
\nabla \left( s(z) \right)_{j+1} &= (j+1) \left( s(z) \right)_j,\quad j \in \mathbb{N},\notag
\end{align}
where $\Pochh{s(z)}{j+1}=s(z)(s(z)+1)\ldots(s(z)+j)$, $\left( s(z) \right)_0=1$, denotes the Pochhammer symbol. In \eqref{relation_delta_nabla}-\eqref{nth_nabla}, when $s(z)=z$, the operator $\nabla$  coincides with the operator $\bigtriangledown$ given in \eqref{triang_nabla}. In what follows we take $s(z)=z^r$.

For type II Charlier-Angelesco polynomial $C_{\vec{n}}^{\vec{a}}(z)$ with multi-index
$\vec{n}\in\mathbb{Z}_{+}^{r}$ and degree
$\lvert\vec{n}\rvert=n_0+\cdots+n_{r-1}$ in $z^r$, the orthogonality relations \eqref{eq:discrete_orth_conditions} on the $r$-star \eqref{r-star} are:

\begin{align}\label{eq:ch_orth_cond_2}
\displaystyle\sum_{k=0}^{\infty} C^{\vec{a}}_{\vec{n}}(z_{\ell,k})
\Pochh{-z_{\ell,k}^r}{j} \rho_{\ell}(z_{\ell,k}) &= 0,\quad 0\leq j
\leq n_{\ell}-1,\quad 0\leq \ell  \leq r-1,
\end{align}
where $\rho_{\ell} (z) = \dfrac{a_{\ell}^{z^{r}}}{\Gamma(z^{r}+1)}$ and $-\pi<\arg a_{\ell}\leq \pi$. 

For type II Meixner-Angelesco polynomials of the second kind $M_{\vec{n}}^{\beta,\vec{c}}(z)$with multi-index
$\vec{n}\in\mathbb{Z}_{+}^{r}$ and degree
$\lvert\vec{n}\rvert=n_0+\cdots+n_{r-1}$ in $z^r$, the orthogonality relations \eqref{eq:discrete_orth_conditions} on the $r$-star \eqref{r-star} are:
\begin{align}\label{eq:m2_orth_cond_2}
\displaystyle\sum_{k=0}^{\infty} M_{\vec{n}}^{\beta,\vec{c}}(z_{\ell,k}) \Pochh{-z_{\ell,k}^r}{j} \rho_{\ell}^{\beta}(z_{\ell,k}) &= 0,\quad 0\leq j  \leq n_{\ell}-1,\quad 0\leq \ell  \leq r-1,
\end{align}
where $\rho_{\ell}^\beta (z) =\dfrac{\Gamma(z^r+\beta)c_{\ell}^{z^{r}}}{\Gamma(z^{r}+1)}$, $-\pi<\arg c_{\ell}\leq \pi$, $\left\lvert c_{\ell} \right\rvert\leq c < 1$, and $(z^{r}_{j,k}+\beta)\in\mathbb{C}\setminus\mathbb{Z}_{-}$.

Using summation by parts in the orthogonality relations \eqref{eq:ch_orth_cond_2} and \eqref{eq:m2_orth_cond_2} one can find the raising operations
\begin{align*} 
&\dfrac{-a_{\ell}}{\rho_{\ell}(z)}\nabla\left( \rho_{\ell}(z) C^{\vec{a}}_{\vec{n}}(z) \right)=C^{\vec{a}}_{\vec{n}+\vec{e}_{\ell}} (z),\\
&
\left( \dfrac{c_{\ell}}{c_{\ell}-1} \right) \dfrac{1}{\rho_{\ell}^{\beta-1}(z)}\nabla\left( \rho_{\ell}^\beta(z) M^{\beta,\vec{c}}_{\vec{n}} (z) \right)=M^{\beta-1,\vec{c}}_{\vec{n}+\vec{e}_{\ell}} \left(z\right). 
\end{align*}

A repeated application of these raising operators gives the following Rodrigues-type formulas
\begin{align*}
C^{\vec{a}}_{\vec{n}} (z) &= \left(\prod_{\ell=0}^{r-1}
\left(-a_{\ell}\right)^{n_{\ell}}\right) \Gamma(z^r+1)
\prod_{j=0}^{r-1} \left( \dfrac{1}{a_j^{z^r}}\nabla^{n_j}a_j^{z^r}
\right) \left(\dfrac{1}{\Gamma(z^r+1)}\right),\\
M^{\beta, \vec{c}}_{\vec{n}} (z)&= \prod_{\ell=0}^{r-1} \left(\dfrac{c_{\ell}}{c_{\ell}-1}\right)^{n_\ell} \dfrac{\Gamma(z^r+1) }{\Gamma(z^r+\beta) }\prod_{j=0}^{r-1} \left( \dfrac{1}{c_j^{z^r}}\nabla^{n_j}c_j^{z^r} \right)\dfrac{\Gamma(z^r+\beta + \lvert\vec{n}\rvert)}{\Gamma(z^r+1)}.
\end{align*}

\section{Type II Meixner-Angelesco polynomials of the second kind}\label{sec:meix2_examples}

Now, we introduce the system of measures for the Meixner-Angelesco polynomials of the second kind. We will find the raising operators, Rodrigues-type formula, and explicit expression for these polynomials. Finally, a different proof the above Theorem \ref{theorem:recurrence_relation} and the coefficients for the nearest neighbor recurrence relations will be given.

Let $c\in \mathbb{C}\setminus\{0\}$, the mass points
\begin{equation} \label{eq:mass_points}
\{z_{j,k}\,\,:\,\,z_{j,k} = k^{1/r}\omega^j\}_{k=0}^\infty \subset \Omega_j,\quad j=0,\ldots,r-1,
\end{equation}
and the weight function,
\begin{equation}\label{eq:ch_weight}
\rho (z) = \dfrac{c^{z^{r}}}{\Gamma(z^{r}+1)}=\dfrac{e^{z^{r}\ln c}}{\Gamma(z^{r}+1)},
\end{equation}  
where $-\pi<\arg c\leq \pi$ (principal branch of the logarithmic function). This function $\rho (z)$ is an extension of the Poisson distribution $\rho(k)=c^{k}/k!$, $c>0$, on the non-negative integers $k=0,1,2\ldots$ (see Meixner case in Table \ref{table1}).

Consider the function $\rho(z)$ in \eqref{eq:ch_weight} and the mass points \eqref{eq:mass_points}, and the complex weight functions
\begin{equation}\label{eq:m1_weight_function}
\rho^{\beta_{j}} (z) =\Gamma(z^r+\beta_{j})\rho (z),\quad -\pi<\arg c\leq \pi,\quad j=0,\ldots, r-1,
\end{equation}  
where the complex parameters $\beta_{j}$ are all different and such that $(z^{r}_{j,k}+\beta_{j})\in\mathbb{C}\setminus\mathbb{Z}_{-}$ and $\beta_i - \beta_j \notin \mathbb{Z}$ for $i \neq j$. In the expression \eqref{eq:m1_weight_function} a complex replacement of the real parameters involved in the negative binomial distribution (Pascal distribution) $\rho^{\beta}(k)=(\beta)_{k}c^{k}/k!$, on the non-negative integers $k=0,1,2\ldots$ is carried out (see Tables \ref{tab1} as well as \cite{bib:nikiforov-uvarov-suslov} and \cite{bib:arvesu-coussement-vanassche} for the real weight functions involving the Meixner polynomials and multiple Meixner polynomials, respectively).

From \eqref{eq:particular_measures}, \eqref{eq:mass_points}, and \eqref{eq:m1_weight_function} we obtain a particular system of discrete measures of Meixner of the second kind $\{ \mu^{\beta_{0}}, \ldots, \mu^{\beta_{r-1}}\}$ supported on the $r$-star \eqref{r-star} with the following moments \eqref{eq:particular_moments}
\begin{equation}
m_{j}^{(\ell)} = \sum_{k=0}^\infty a_{k,j}^{(\ell)},\quad
a_{k,j}^{(\ell)}
=k^j \dfrac{\Gamma(k + \beta_{\ell})}{\Gamma(k + 1)}c^k,\quad \beta_{\ell}\not=0,-1,-2,\ldots,\quad \ell=0,\ldots,r-1,\quad j\in \mathbb{N}.
\label{eq:m1_moments}
\end{equation}
For any fixed $j \in \mathbb{N}$ and $\ell = 0,...,r-1$, one gets 
\begin{equation*}
\left\lvert
a_{k+1,j}^{(\ell)}/a_{k,j}^{(\ell)}\right\rvert = \left\lvert c \right\rvert \left\lvert \dfrac{(k+\beta_{\ell})(k+1)^{j-1}}{k^j} \right\rvert \longrightarrow \left\lvert c \right\rvert, \ \ \text{as} \ \ k \rightarrow \infty.
\end{equation*} 

Taking $\left\lvert c \right\rvert\leq a < 1$ for all $\ell = 0, \ldots, r-1$, from the Ratio Test follows the absolute convergence of the series. Then, \eqref{eq:m1_moments} converges and all moments exist. Moreover, the matrix $\mathcal{M}$ derived from the linear system \eqref{eq:discrete_orth_conditions} has full rank (see expressions \eqref{matrix} and \eqref{Pn_matrix_eq}). Thus, the system of discrete measures of Meixner-Angelesco of the second kind is a perfect system of measures.

\begin{definition} The type II monic Meixner-Angelesco polynomial of the second kind $M_{\vec{n}}^{\vec{\beta},c}$ for the multi-index $\vec{n}=(n_0,\ldots,n_{r-1})$ 
on the $r$-star and weight functions \eqref{eq:m1_weight_function} is the polynomial of degree 
$\lvert\vec{n}\rvert$ in $z^r$ defined by the orthogonality relations

\begin{equation} \label{eq:m1_orth_cond_1}
\displaystyle{\int\limits_{\Omega_\ell}} M_{\vec{n}}^{\vec{\beta},c}(z) \left( -z^r \right)_j d\mu^{\beta_{\ell}}=0, \quad j=0,\ldots, n_{\ell}-1,\quad \ell=0,\ldots,r-1,
\end{equation}
or equivalently,
\begin{align}\label{eq:m1_orth_cond_2}
\displaystyle\sum_{k=0}^{\infty} M_{\vec{n}}^{\vec{\beta},c}(z_{\ell,k}) \Pochh{-z_{\ell,k}^r}{j} \rho^{\beta_{\ell}}(z_{\ell,k}) &= 0,\quad 0\leq j  \leq n_{\ell}-1,\quad 0\leq \ell  \leq r-1.
\end{align}
\end{definition}

The orthogonality relations \eqref{eq:m1_orth_cond_1} and \eqref{eq:m1_orth_cond_2} are a particular situation of \eqref{eq:discrete_orth_conditions}, with the $\omega$-symmetric weight functions \eqref{eq:m1_weight_function}. The system \eqref{eq:m1_orth_cond_2} defines $\lvert\vec{n}\rvert$ conditions for the $\lvert\vec{n}\rvert$-unknown coefficients of the monic polynomial $M_{\vec{n}}^{\vec{\beta},c}$ of degree $\lvert\vec{n}\rvert$ in $z^r$.

\subsection{Raising relation and Rodrigues-type formula}\label{raising}

Here we study the \textit{raising operator} and \textit{Rodrigues-type formula} for the type II Meixner-Angelesco polynomials of the second kind.

Define 
\begin{equation}\label{Psi_operator}
\Psi^{\beta_\ell} =\left( \dfrac{c}{c-1} \right) \dfrac{1}{\rho^{\beta_{\ell}-1}(z)}\nabla \rho^{\beta_{\ell}}(z),
\end{equation}
where $\rho^{\beta_{\ell}}(z)$ is given in \eqref{eq:m1_weight_function}. For any power $k=0,1,\ldots$, 
\begin{align}\label{Psi_power_m1}
\Psi^{\beta_{\ell}} z^k&=\dfrac{c}{c-1} \dfrac{c^{-z^r}\Gamma(z^r+1)}{\Gamma(z^r+\beta_{\ell}-1)}
\left(\nabla \dfrac{c^{z^r}\Gamma(z^r+\beta_{\ell})}{\Gamma(z^r+1)}z^k\right)\notag\\
&=\dfrac{c}{c-1} \dfrac{\Gamma(z^r+1)}{\Gamma(z^r+\beta_{\ell}-1)}\left(\dfrac{\Gamma(z^r+\beta_{\ell})z^k}{\Gamma(z^r+1)}-
\dfrac{c^{-1}(z^k-1+\delta_{k,0})\Gamma(z^r+\beta_{\ell}-1)}{\Gamma(z^r)}\right)\notag\\
&=z^{r+k}+\dfrac{c(\beta_{\ell}-1)}{c-1}z^{k}+\dfrac{(1-\delta_{k,0})}{c-1}z^r,
\end{align}
where $\delta_{k,0}$ denotes the Kronecker delta. Notice that $\deg\Psi^{\beta_{\ell}}z^k=r+k$ and for $k=r\lvert\vec{n}\rvert$ the resulting polynomial has degree $r(\lvert\vec{n}\rvert+1)$. For polynomials $p(z)\in\mathbb{P}_{N}=\spn\{1,z,\ldots,z^{N}\}$, we will use the expression \eqref{Psi_power_m1} termwise.

\begin{lemma} \label{lemma:L1_conmutativity}
The following commutative property 
\begin{equation}\label{operator_comm_L1}
\Psi^{\beta_{j}}\Psi^{\beta_{\ell}}=\Psi^{\beta_{\ell}}\Psi^{\beta_{j}},\quad j\not=\ell,\quad 0\leq j,\ell\leq r-1,
\end{equation}
holds.
\end{lemma}

\begin{proof} For any power $k\in\mathbb{Z}_{+}$, using \eqref{Psi_power_m1} first in the index $\ell$ and then in $j$, one gets
\begin{align*}
\Psi^{\beta_{j}}\Psi^{\beta_{\ell}}z^k&=\Psi^{\beta_{j}}z^{r+k}+\dfrac{c(\beta_{\ell}-1)}{c-1}\Psi^{\beta_{j}}z^{k}+\dfrac{1}{c-1}\Psi^{\beta_{j}}z^r\\
&=z^{2r+k}+\dfrac{c}{c-1}\left[(\beta_{j}-1)+(\beta_{\ell}-1)\right]z^{r+k}+\dfrac{c^2}{(c-1)^2}(\beta_{j}-1)(\beta_{\ell}-1)z^{k}\\&
+\dfrac{(1-\delta_{k,0})}{c-1}z^{2r}+\dfrac{\left(c(\beta_{\ell}+\beta_{j})-(c+\delta_{k,0})\right)(1-\delta_{k,0})}{(c-1)^2}z^{r}.
\end{align*}
Clearly, the right-hand side of this equation is invariant under change of complex parameters $\beta_{\ell}\leftrightarrow\beta_{j}$, which implies \eqref{operator_comm_L1}.
\end{proof}

\begin{lemma} \label{lemma:m1_raising_operator}
The type II Meixner-Angelesco polynomials of the second kind satisfy the raising relation
\begin{equation} \label{eq:m1_raising}
M^{\vec{\beta}-\vec{e}_{\ell},c}_{\vec{n}+\vec{e}_{\ell}} \left(z\right) = \Psi^{\beta_{\ell}} M^{\vec{\beta},c}_{\vec{n}} \left(z \right),
\end{equation}
where $\Psi^{\beta_{\ell}}$ is given in \eqref{Psi_operator}.
\end{lemma}

The operator $\Psi^{\beta_{\ell}}$ is called \textit{raising operator} because the ${\ell}$-th component of the multi-index $\vec{n}$ is increased by $1$.

\begin{proof} Replacing $(-z_{{\ell},k}^r)_{j}$ by $\nabla(-z_{{\ell},k}^r +1)_{j+1}$, $j=0,\ldots,n_{\ell}-1$, in the orthogonality relation \eqref{eq:m1_orth_cond_2} one has
\begin{equation}\label{orth_transition_m1}
\displaystyle{\sum^\infty_{k=0}} M^{\vec{\beta},c}_{\vec{n}}\left(z_{{\ell},k}\right) \nabla\left(-z_{{\ell},k}^r +1\right)_{j+1} \rho^{\beta_{\ell}} \left(z_{{\ell},k} \right)=0,\quad 0\leq \ell  \leq r-1.
\end{equation}
Note that $\nabla(-z^r+1)_{j+1}$ is a polynomial of degree $j$ in $z^r$.

Considering the explicit expression of the mass points \eqref{eq:mass_points} along the rays, one gets $z_{{\ell},k}^r = k$, for ${\ell} = 0, \ldots, r-1$, hence the orthogonality relations \eqref{orth_transition_m1} become 
\begin{equation}\label{orth_new_m1}
\displaystyle\sum^\infty_{k=0} M^{\vec{\beta},c}_{\vec{n}}(k) 
\nabla\left(-k +1\right)_{j+1} \rho^{\beta_{\ell}}_{k}=0,\quad 0\leq j  \leq n_{\ell}-1,\quad 0\leq \ell  \leq r-1,
\end{equation}
where $\rho^{\beta_{\ell}}_{k}=c^{k}\Gamma(\beta_{\ell}+k)/\Gamma(k+1)$.

In \eqref{orth_new_m1} we will use summation by parts and the relations
\begin{align}
\rho^{\beta_{\ell}}_{-1} = c^{-1}\Gamma(\beta_{\ell}-1)/\Gamma\left(0\right) &=0,\notag\\
\lim_{k \rightarrow \infty} M^{\vec{\beta},c}_{\vec{n}}\left(k \right) \left( - k  \right)_{j+1} \rho^{\beta_{\ell}}_{k}&=0,\label{lim_moment_term_m1}
\end{align}
where \eqref{lim_moment_term_m1} follows from the $n$th-Term Test for the absolutely convergent series \eqref{eq:m1_moments}. Moreover, from \eqref{relation_delta_nabla} the relations \eqref{orth_new_m1} become
\begin{align*}
0=\sum^\infty_{k=0} M^{\vec{\beta},c}_{\vec{n}}(k) 
\nabla\left(-k +1\right)_{j+1} \rho^{\beta_{\ell}}_{k}
&=\sum_{k=0}^{\infty}
M^{\vec{\beta},c}_{\vec{n}}(k)\rho^{\beta_{\ell}}_{k}\Delta \left(-k\right)_{j+1}
\\
&=-\sum_{k=0}^{\infty}\left(-k\right)_{j+1}
\nabla\left(M^{\vec{\beta},c}_{\vec{n}}(k)
\rho^{\beta_{\ell}}_{k}\right).
\end{align*}
Note that from \eqref{Psi_operator}-\eqref{Psi_power_m1}
\begin{equation}\label{eq4recurrence_m1}
\nabla\left(M^{\vec{\beta},c}_{\vec{n}}(k)
\rho^{\beta_{\ell}}_{k}\right)=\mathcal{Q}_{\lvert\vec{n}\rvert+1}(k)\rho^{\beta_{\ell}-1}_{k},
\end{equation}
where $\mathcal{Q}_{\lvert\vec{n}\rvert+1}(k)=(k+\beta_{\ell}-1) M^{\vec{\beta},c}_{\vec{n}}  \left(k \right) - \dfrac{k}{c}   M^{\vec{\beta},c}_{\vec{n}} \left(k -1\right)$, which is a polynomial of degree exactly $\lvert\vec{n}\rvert+1$. Therefore, from the uniqueness of the monic multiple orthogonal polynomials derived from the orthogonality relations \eqref{orth_transition_m1} one gets $\mathcal{Q}_{\lvert\vec{n}\rvert+1}(z)=\left(\dfrac{c-1}{c}\right) M^{\vec{\beta}-\vec{e}_{\ell},c}_{\vec{n}+\vec{e}_{\ell}} \left( z \right)$, which is a polynomial of degree $\lvert\vec{n}\rvert+1$ in $z^r$. The expression \eqref{eq4recurrence_m1} becomes
$$
\nabla\left(M^{\vec{\beta},c}_{\vec{n}}(z)
\rho^{\beta_{\ell}}(z)\right)=\left(\dfrac{c-1}{c}\right) M^{\vec{\beta}-\vec{e}_{\ell},c}_{\vec{n}+\vec{e}_{\ell}} (z)\rho^{\beta_{\ell}-1}(z).
$$
Finally, taking into account \eqref{Psi_operator} the raising relation \eqref{eq:m1_raising} holds.
\end{proof}

\begin{remark} Observe that from \eqref{eq:m1_raising} and \eqref{eq4recurrence_m1} one has the following recurrence relation involving the type II Meixner-Angelesco polynomials in $z^r$,
$$
M^{\vec{\beta}-\vec{e}_{\ell},c}_{\vec{n}+\vec{e}_{\ell}} (z)=\left( \dfrac{c}{c-1} \right)\left((z^r+\beta_{\ell}-1) M^{\vec{\beta},c}_{\vec{n}}  \left(z \right) - \dfrac{1}{c} z^r  M^{\vec{\beta},c}_{\vec{n}} \left(z -1\right)\right).  
$$
\end{remark}

In the sequel we will deal with the compositions of raising operators involving the expression given in \eqref{Psi_operator}. Thus, define 
\begin{align}\label{eq:m1_kth_times_operator_Psi}
\mathcal{R}^{\beta_{\ell}}_{k}&=\Psi^{\beta_{\ell}+1}\circ\cdots\circ \Psi^{\beta_{\ell}+k}\notag\\
&= \left( \dfrac{c}{c-1} \right)^{k}\left( \dfrac{1}{\rho^{\beta_{\ell}}(z)}\nabla^{k} \rho^{\beta_{\ell}+k} (z)\right)\notag\\
&=\left( \dfrac{c}{c-1} \right)^{k}\dfrac{\Gamma(z^r+1) }{c^{z^r} }\left( \dfrac{1}{\Gamma(z^r+\beta_{\ell})}\nabla^{k} \Gamma(z^r+\beta_{\ell}+k )\right)\left(\dfrac{c^{z^r}}{\Gamma(z^r+1)}\right).
\end{align}
Moreover, from Lemma \ref{lemma:L1_conmutativity} the following commutative relation
\begin{equation}\label{commutation_R_m-k}
\mathcal{R}^{\beta_{\ell}}_{k}\mathcal{R}^{\beta_{j}}_{m}=\mathcal{R}^{\beta_{j}}_{m}\mathcal{R}^{\beta_{\ell}}_{k},\qquad k,m\in\mathbb{Z}_{+},\quad 0\leq j,\ell\leq r-1,
\end{equation}
holds.
\begin{theorem}\label{theorem:m1_rodrigues} The type II Meixner-Angelesco polynomials of the second kind verify the Rodrigues-type formula

\begin{align}
M^{\vec{\beta},c}_{\vec{n}} (z)&= \left(\dfrac{c}{c-1}\right)^{|\vec{n}|} \dfrac{\Gamma(z^r+1) }{c^{z^r}}\prod_{j=0}^{r-1} 
\left(\dfrac{1}{\Gamma(z^r+\beta_j)}\nabla^{n_j} \Gamma(z^r+\beta_j + n_j)\right)\dfrac{c^{z^r}}{\Gamma(z^r+1)} \label{eq:m1_rodrigues_util}\\ 
&= \left(\dfrac{c}{c-1}\right)^{|\vec{n}|}\prod_{j=0}^{r-1} \left( \dfrac{1}{\rho^{\beta_j}(z)}\nabla^{n_j} \rho^{\beta_j + n_j}(z) \right).\notag
\end{align}
\end{theorem}

\begin{proof} In \eqref{eq:m1_raising} replace $\beta_{\ell}$ by $\beta_{\ell}+n_{\ell}$ (for any fixed $\ell$-ray, $\ell\in\{0,\ldots, r-1\}$). Hence,
\begin{align}\label{m1_initial_exp_rod}
M^{\vec{\beta}+(n_{\ell}-1)\vec{e}_{\ell}, c}_{\vec{e}_{\ell}} (z) &= \Psi^{\beta_{\ell}+n_{\ell}} M^{\vec{\beta}+n_{\ell}\vec{e}_{\ell}, c}_{\vec{0}} (z)= \left(\dfrac{c}{c-1} \right)\left(\dfrac{1}{\rho^{\beta_{\ell}+n_{\ell}-1}(z)}\nabla \rho^{\beta_{\ell}+n_{\ell}}(z)\right),
\end{align}
where $M^{\vec{\beta}+n_{\ell}\vec{e}_{\ell}, c}_{\vec{0}} (z)=1$

Multiplying equation \eqref{m1_initial_exp_rod} from the left by the product of $(n_{\ell}-1)$-raising operators $\Psi^{\beta_{\ell}+1}\cdots \Psi^{\beta_{\ell}+n_{\ell}-1}$ and using \eqref{eq:m1_kth_times_operator_Psi}, one gets
\begin{align}\label{eq:m1_rodrigues_proof_one_ray}
M^{\vec{\beta},c}_{n_{\ell}\vec{e}_{\ell}} (z)&=\prod\limits_{t=0}^{n_{\ell}-1} \Psi^{\beta_{\ell}+1 + t}=\left(\dfrac{c}{c-1}\right)^{n_{\ell}} \left( \dfrac{1}{\rho^{\beta_{\ell}}(z)}\nabla^{n_{\ell}} \rho^{\beta_{\ell}+n_{\ell}}(z)\right). 
\end{align}
Here we have used Lemma \ref{lemma:m1_raising_operator} and the orthogonality conditions \eqref{eq:m1_orth_cond_2} for the multi-index $n_{\ell}\vec{e}_\ell$. 

Now, take a ray $j\not=\ell$, $j=0,\ldots,r-1$, and replace in \eqref{eq:m1_rodrigues_proof_one_ray} $\beta_j$ by $\beta_j+n_j$, then multiply \eqref{eq:m1_rodrigues_proof_one_ray} from the left by the product of $n_{j}$-raising operators $\Psi^{\beta_{j}+1}\cdots \Psi^{\beta_{j}+n_{j}}$ to obtain
\begin{multline}
M^{\vec{\beta},c}_{n_{j}\vec{e}_{j}+n_{\ell}\vec{e}_{\ell}} (z)
=\left(\prod\limits_{s=0}^{n_{j}-1} \Psi^{\beta_{j}+1 + s}\right) M^{\vec{\beta}+n_{j}\vec{e}_{j}, c}_{n_{\ell}\vec{e}_{\ell}} (z)\label{eq:m1_rodrigues_proof_second_ray}\\
=\left(\prod\limits_{s=0}^{n_{j}-1} \Psi^{\beta_{j}+1 + s}\right)
\left(\prod\limits_{t=0}^{n_{\ell}-1} \Psi^{\beta_{\ell}+1 + t}\right)M^{\vec{\beta}+n_{j}\vec{e}_{j}+n_{\ell}\vec{e}_{\ell}, c}_{\vec{0}} (z)\\
=\dfrac{ c^{n_{j}+n_{\ell}} }{ (c-1)^{n_{j}+n_{\ell}} }
\left( \dfrac{1}{\rho^{\beta_{j}}(z)}\nabla^{n_{j}} \rho^{\beta_{j}+n_{j}}(z)\right)
 \left( \dfrac{1}{\rho^{\beta_{\ell}}(z)}\nabla^{n_{\ell}} \rho^{\beta_{\ell}+n_{\ell}}(z)\right), 
\end{multline}
where $M^{\vec{\beta}+n_{j}\vec{e}_{j}+n_{\ell}\vec{e}_{\ell}, c}_{\vec{0}} (z)=1$. Hence, using the explicit expression of the weight functions in \eqref{eq:m1_weight_function} one has
\begin{align}\label{eq:m1_rodrigues_proof_second_ray_equiv}
M^{\vec{\beta},c}_{n_{j}\vec{e}_{j}+n_{\ell}\vec{e}_{\ell}} (z)
&=\left(\dfrac{c}{c-1}\right)^{n_{j}+n_{\ell}}\dfrac{\Gamma(z^r+1) }{c^{z^r}}\left(\dfrac{1}{\Gamma(z^r+\beta_j)}\nabla^{n_j} \Gamma(z^r+\beta_j + n_j)\right)\notag\\
&\times\left(\dfrac{1}{\Gamma(z^r+\beta_{\ell})}\nabla^{n_{\ell}} \Gamma(z^r+\beta_{\ell} + n_{\ell})\right)
\dfrac{c^{z^r}}{\Gamma(z^r+1)}.
\end{align}

From Lemma \ref{lemma:L1_conmutativity} the above raising operators are commuting so they can be taken in any order. One can begin with the $j$-ray, $j\in\{0,\ldots, r-1\}$ and then continue with other ray $\ell\not=j$, $\ell\in\{0,\ldots, r-1\}$, leading to the same expressions \eqref{eq:m1_rodrigues_proof_second_ray}-\eqref{eq:m1_rodrigues_proof_second_ray_equiv}. By continuing this process with all other rays taken in any order one gets \eqref{eq:m1_rodrigues_util}. Indeed, if $\kappa=j_0,\ldots,j_{r-1}$ is one of the $r!$ permutations of $\{0,\ldots,r-1\}$ and $\vec{n}_{\kappa}=(n_{j_0},\ldots,n_{j_{r-1}})$, the corresponding multi-index from $\vec{n}=(n_0,\ldots,n_{r-1})$, where $\lvert\vec{n}_{\kappa}\rvert=\lvert\vec{n}\rvert$, one has
\begin{multline}\label{eq:m1_4_theorem_rodrigues}
M^{\vec{\beta}_{\kappa},c}_{\vec{n}_{\kappa}} (z)=\prod\limits_{t_{0}=0}^{n_{j_{0}}-1} \Psi^{\beta_{j_{0}}+1 + t_{0}}
\prod\limits_{t_{1}=0}^{n_{j_{1}}-1} \Psi^{\beta_{j_{1}}+1 + t_{1}}\cdots
\prod\limits_{t_{r-1}=0}^{n_{j_{r-1}}-1} \Psi^{\beta_{j_{r-1}}+1 + t_{r-1}}=\left(\dfrac{c}{c-1}\right)^{|\vec{n}|} \dfrac{\Gamma(z^r+1) }{c^{z^r}}\\
\times\left(\dfrac{1}{\Gamma(z^r+\beta_{j_{0}})}\nabla^{n_{j_{0}}} \Gamma(z^r+\beta_{j_{0}} + n_{j_{0}})\right)
\cdots\left(\dfrac{1}{\Gamma(z^r+\beta_{j_{r-1}})}\nabla^{n_{j_{r-1}}} \Gamma(z^r+\beta_{j_{r-1}} + n_j)\right)\dfrac{c^{z^r}}{\Gamma(z^r+1)}.
\end{multline}
Here $\vec{\beta}_{\kappa}=(\beta_{j_0},\ldots,\beta_{j_{r-1}})$ is the permutation from $\{\beta_1,\ldots,\beta_{r-1}\}$ associated with the above multi-index $\vec{n}_{\kappa}=(n_{j_0},\ldots,n_{j_{r-1}})$.

From the commutative relation \eqref{operator_comm_L1} in Lemma \ref{lemma:L1_conmutativity}, the above product in \eqref{eq:m1_4_theorem_rodrigues}
$$
\left(\dfrac{1}{\Gamma(z^r+\beta_{j_{0}})}\nabla^{n_{j_{0}}} \Gamma(z^r+\beta_{j_{0}} + n_{j_{0}})\right)
\cdots\left(\dfrac{1}{\Gamma(z^r+\beta_{j_{r-1}})}\nabla^{n_{j_{r-1}}} \Gamma(z^r+\beta_{j_{r-1}} + n_j)\right),
$$ 
can be expressed (reordered) as follows
$$
\prod_{j=0}^{r-1} 
\left(\dfrac{1}{\Gamma(z^r+\beta_j)}\nabla^{n_j} \Gamma(z^r+\beta_j + n_j)\right).
$$
Therefore, the following relation holds
\begin{align*}
M^{\vec{\beta}_{\kappa},c}_{\vec{n}_{\kappa}} (z)&=\prod\limits_{t_{0}=0}^{n_{j_{0}}-1} \Psi^{\beta_{j_{0}}+1 + t_{0}}
\prod\limits_{t_{1}=0}^{n_{j_{1}}-1} \Psi^{\beta_{j_{1}}+1 + t_{1}}\cdots
\prod\limits_{t_{r-1}=0}^{n_{j_{r-1}}-1} \Psi^{\beta_{j_{r-1}}+1 + t_{r-1}} \\
&=\prod\limits_{t_{0}=0}^{n_{0}-1} \Psi^{\beta_{0}+1 + t_{0}}
\prod\limits_{t_{1}=0}^{n_{1}-1} \Psi^{\beta_{1}+1 + t_{1}}\cdots
\prod\limits_{t_{r-1}=0}^{n_{r-1}-1} \Psi^{\beta_{r-1}+1 + t_{r-1}}\\
&=\left(\dfrac{c}{c-1}\right)^{|\vec{n}|} \dfrac{\Gamma(z^r+1) }{c^{z^r}}\prod_{j=0}^{r-1} 
\left(\dfrac{1}{\Gamma(z^r+\beta_j)}\nabla^{n_j} \Gamma(z^r+\beta_j + n_j)\right)\dfrac{c^{z^r}}{\Gamma(z^r+1)}\\
&=\left(\dfrac{c}{c-1}\right)^{|\vec{n}|}\prod_{j=0}^{r-1} \left( \dfrac{1}{\rho^{\beta_j}(z)}\nabla^{n_j} \rho^{\beta_j + n_j}(z) \right)\\
&=M^{\vec{\beta},c}_{\vec{n}} (z).
\end{align*}
This concludes the proof of \eqref{eq:m1_rodrigues_util}.
\end{proof}

\begin{remark}
The formula \eqref{eq:m1_rodrigues_proof_one_ray} gives the explicit expression for the Meixner polynomials of degree $n_{\ell}\in\mathbb{Z}_{+}$, in the variable $z^r$, which are orthogonal with respect to the discrete measure derived from the complex weight $\rho^{\beta_{\ell}}(z)$ on the $\ell$-ray, for $\ell\in\{0,\ldots, r-1\}$. In this situation $\vec{n}=(0,\ldots,n_{\ell},\ldots,0)$, that is, one is dealing with scalar orthogonality. Notice that the expression \eqref{eq:m1_rodrigues_proof_one_ray} is similar to \eqref{meixner_classic} for the classical Meixner polynomials. The similarity is understood as a complex replacement of the real parameters involved in the orthogonality measure as well as the change of the real variable $x$ by $z^r$ (see formula \eqref{eq:m1_explicit_expression}).  
\end{remark}

\section{Explicit expression and recurrence relation}\label{recu}

Here we obtain the explicit expression for the type II Meixner-Angelesco polynomials of the second kind from the Rodrigues formula \eqref{eq:m1_rodrigues_util}. As a consequence of Theorem \ref{theorem:m1_rodrigues} we have the following Corollary.

\begin{corollary}\label{m1_quasi_final_corollary} The type II Meixner-Angelesco polynomials of the second kind on the $r$-star are given by 
\begin{equation} \label{eq:m1_explicit_expression}
M^{\vec{\beta},c}_{\vec{n}}(z) = \dfrac{c^{|\vec{n}|}}{\left( c-1 \right)^{|\vec{n}|}}
\sum_{k_0=0}^{n_{0}}\cdots \sum_{k_{r-1}=0}^{n_{r-1}}
{n_{0}\choose{k_{0}}}\cdots {n_{\ell}\choose{k_{\ell}}}
\dfrac{\Pochh{-z^r}{\lvert\vec{k}\rvert}}{c^{\lvert\vec{k}\rvert}} 
\prod_{m=0}^{r-1}\Pochh{z^r+\beta_{m}-\lvert\vec{k}\rvert_{m}+k_{m}}{n_{m}-k_{m}},
\end{equation}
where $\lvert\vec{k}\rvert = k_0 + \ldots + k_{r-1}$ and $\lvert\vec{k}\rvert_{m}=\sum_{j=0}^{m}k_{j}$.
\end{corollary} 

\begin{proof} Using formula \eqref{nth_nabla} in \eqref{eq:m1_kth_times_operator_Psi} yields
\begin{multline}\label{particular_operator_R}
\mathcal{R}_{n_{\ell}}^{\beta_{\ell}}
=\left(\dfrac{c}{c-1}\right)^{n_{\ell}}\dfrac{\Gamma(z^r+1) }{c^{z^r}}\dfrac{1}{\Gamma(z^r+\beta_{\ell})}\nabla^{n_{\ell}} \Gamma(z^r+\beta_{\ell} + n_{\ell})\dfrac{c^{z^r}}{\Gamma(z^r+1)}\\
=\left(\dfrac{c}{c-1}\right)^{n_{\ell}}\dfrac{\Gamma(z^r+1) }{c^{z^r}}\sum_{k_{\ell}=0}^{n_{\ell}} {n_{\ell}\choose{k_{\ell}} } (-1)^{k_{\ell}}c^{z^r-k_{\ell}}\dfrac{\Gamma(z^r+\beta_{\ell}+n_{\ell}-k_{\ell})}{\Gamma(z^r+1-k_{\ell})\Gamma(z^r+\beta_{\ell})}\\
=\left(\dfrac{c}{c-1}\right)^{n_{\ell}}\dfrac{\Gamma(z^r+1) }{c^{z^r}}\sum_{k_{\ell}=0}^{n_{\ell}} {n_{\ell}\choose{k_{\ell}} } (-1)^{k_{\ell}}c^{z^r-k_{\ell}}\dfrac{(z^r+\beta_{\ell})_{n_{\ell}-k_{\ell}}}{\Gamma(z^r+1-k_{\ell})},
\end{multline}
where the following relation has been used
\begin{align*}
\dfrac{\Gamma(z^r+\beta_{\ell}+n_{\ell}-k_{\ell})}{\Gamma(z^r+\beta_{\ell})}
&=\Pochh{z^r+\beta_{\ell}}{n_{\ell}-k_{\ell}}.
\end{align*}
Moreover, taking into account the relation $
\dfrac{\Gamma(z^r+1)}{\Gamma(z^r+1-k_{\ell})}=(-1)^{k_{\ell}}
\Pochh{-z^r}{k_{\ell}}$, equation \eqref{particular_operator_R} becomes
\begin{equation*}\label{R_nl_explicit}
\mathcal{R}_{n_{\ell}}^{\beta_{\ell}}
=\left(\dfrac{c}{c-1}\right)^{n_{\ell}}\sum_{k_{\ell}=0}^{n_{\ell}} {n_{\ell}\choose{k_{\ell}} } c^{-k_{\ell}}\Pochh{-z^r}{k_{\ell}}(z^r+\beta_{\ell})_{n_{\ell}-k_{\ell}},
\end{equation*}
which is an explicit expression for the polynomial $M^{\vec{\beta},c}_{n_{\ell}\vec{e}_{\ell}} (z)$ in equation \eqref{eq:m1_rodrigues_proof_one_ray}.

Now, using \eqref{nth_nabla} one has
\begin{align}\label{polymomial-2-indices}
\mathcal{R}_{n_{j}}^{\beta_{j}}\mathcal{R}_{n_{\ell}}^{\beta_{\ell}}&=\left(\dfrac{c}{c-1}\right)^{n_{\ell}}\mathcal{R}_{n_{j}}^{\beta_{j}}\left(\dfrac{\Gamma(z^r+1) }{c^{z^r}}\sum_{k_{\ell}=0}^{n_{\ell}} {n_{\ell}\choose{k_{\ell}} } (-1)^{k_{\ell}}c^{z^r-k_{\ell}}\dfrac{(z^r+\beta_{\ell})_{n_{\ell}-k_{\ell}}}{\Gamma(z^r+1-k_{\ell})}\right)\notag\\
&=\left(\dfrac{c}{c-1}\right)^{n_{j}+n_{\ell}}\dfrac{\Gamma(z^r+1) }{c^{z^r}\,\Gamma(z^r+\beta_{j})}\nabla^{n_{j}} \Gamma(z^r+\beta_{j} + n_{j})
\sum_{k_{\ell}=0}^{n_{\ell}} {n_{\ell}\choose{k_{\ell}} } (-1)^{k_{\ell}}c^{z^r-k_{\ell}}\dfrac{(z^r+\beta_{\ell})_{n_{\ell}-k_{\ell}}}{\Gamma(z^r+1-k_{\ell})}\notag\\
&=\left(\dfrac{c}{c-1}\right)^{n_{j}+n_{\ell}}\dfrac{\Gamma(z^r+1) }{c^{z^r}}
\sum_{k_{j}=0}^{n_{j}} {n_{j}\choose{k_{j}} }(-1)^{k_{j}}(z^r+\beta_{j})_{n_{j}-k_{j}}
\sum_{k_{\ell}=0}^{n_{\ell}} {n_{\ell}\choose{k_{\ell}} } (-1)^{k_{\ell}}c^{z^r-k_{j}-k_{\ell}}\notag\\
&\times\dfrac{(z^r+\beta_{\ell}-k_{j})_{n_{\ell}-k_{\ell}}}{\Gamma(z^r+1-k_{j}-k_{\ell})}\\
&=\left(\dfrac{c}{c-1}\right)^{n_{j}+n_{\ell}}
\sum_{k_{j}=0}^{n_{j}}\sum_{k_{\ell}=0}^{n_{\ell}}{n_{j}\choose{k_{j}} }
{n_{\ell}\choose{k_{\ell}} } c^{-k_{j}-k_{\ell}}\Pochh{-z^r}{k_{j}+k_{\ell}}(z^r+\beta_{j})_{n_{j}-k_{j}}(z^r+\beta_{\ell}-k_{j})_{n_{\ell}-k_{\ell}}\notag.
\end{align}
Notice that this expression coincides with \eqref{eq:m1_rodrigues_proof_second_ray_equiv}.

Finally, repeating the above process, base on the structure of the Rodrigues formula \eqref{eq:m1_rodrigues_util} and Lemma \ref{lemma:L1_conmutativity} one gets
\begin{align*}
M^{\vec{\beta},c}_{\vec{n}}(z) = \prod_{\ell=0}^{r-1}\mathcal{R}_{n_{\ell}}^{\beta_{\ell}}&=\dfrac{c^{|\vec{n}|}}{\left( c-1 \right)^{|\vec{n}|}}
\sum_{k_0=0}^{n_{0}}\cdots \sum_{k_{r-1}=0}^{n_{r-1}}
{n_{0}\choose{k_{0}}}\cdots {n_{\ell}\choose{k_{\ell}}}
\dfrac{\Pochh{-z^r}{\lvert\vec{k}\rvert}}{c^{\lvert\vec{k}\rvert}}\\
&\times\Pochh{z^r+\beta_{0}}{n_0-k_0}\Pochh{z^r+\beta_{1}-k_0}{n_1-k_1}\cdots
\Pochh{z^r+\beta_{r-1}-\lvert\vec{k}\rvert+k_{r-1}}{n_{r-1}-k_{r-1}},
\end{align*}
which coincides with \eqref{eq:m1_explicit_expression}.

\end{proof}

In order to compute the coefficients of the recurrence relation for the type II Meixner-Angelesco polynomials of the second kind given in Theorem \ref{theorem:recurrence_relation} we use the above raising operators. The procedure is similar to the one used in \cite{bib:arvesu-ramirez1}. 

\begin{lemma}\label{auxiliry-lemma} The following relation holds
\begin{equation}\label{eq_lemma_auxi}
\mathcal{R}_{n_{\ell}+1}^{\beta_{\ell}}=z^{r}\mathcal{R}_{n_{\ell}}^{\beta_{\ell}}+ \left(\dfrac{c\left(n_{\ell}+\beta_{\ell}\right)}{c -1}+\dfrac{n_{\ell}}{c -1 }\right)\mathcal{R}_{n_{\ell}}^{\beta_{\ell}}
-\dfrac{c}{\left(c - 1\right)^2} n_{\ell} (\beta_{\ell} + n_{\ell}- 1)\mathcal{R}_{n_{\ell}-1}^{\beta_{\ell}}.
\end{equation}
\end{lemma} 
\begin{proof} We begin by computing the action of the operator $\mathcal{R}_{n_{\ell}}^{\beta_{\ell}}$ on the function $z^{r}$. Use the first relation in \eqref{nth_nabla} $m$-times over the product of functions $\left(z^r\rho^{\beta_{\ell}+k} (z)\right)$, that is, the Leibniz's rule for the {\it n}th derivative to get
\begin{align}\label{pre_eq_R_operator-1}
\nabla^{m} \left(z^r\rho^{\beta_{\ell}+k} (z)\right)=m\nabla^{m-1} \rho^{\beta_{\ell}+k} (z)+(z^r-m)\nabla^{m} \rho^{\beta_{\ell}+k} (z),\quad m\in\mathbb{Z}_{+},
\end{align}
where $\nabla^{0}$ is the identity operator.

Then, in equation \eqref{pre_eq_R_operator-1}, take $k=m=n_{\ell}$ and multiply from the left by $\left( \dfrac{c}{c-1} \right)^{n_{\ell}}\dfrac{1 }{\rho^{\beta_{\ell}} (z)}$ to obtain
\begin{align} 
\mathcal{R}_{n_{\ell}}^{\beta_{\ell}}z^{r}&=(z^r-n_{\ell})\mathcal{R}_{n_{\ell}}^{\beta_{\ell}}+\left( \dfrac{c}{c-1} \right)^{n_{\ell}}\dfrac{n_{\ell}}{\rho^{\beta_{\ell}} (z)}\nabla^{n_{\ell}-1}\rho^{\beta_{\ell}+n_{\ell}} (z)\notag\\
&=\left(z^r+\dfrac{n_{\ell}}{c-1}\right)\mathcal{R}_{n_{\ell}}^{\beta_{\ell}}-\dfrac{c}{\left(c - 1\right)^2} n_{\ell} (\beta_{\ell} + n_{\ell}- 1)\mathcal{R}_{n_{\ell}-1}^{\beta_{\ell}}\label{important_Rn}.
\end{align}
Here the following relation has been used:
\begin{equation}\label{simple_eq_nablas_Rs}
\left( \dfrac{c}{c-1} \right)^{n_{\ell}}\dfrac{n_{\ell}}{\rho^{\beta_{\ell}} (z)}\nabla^{n_{\ell}-1}\rho^{\beta_{\ell}+n_{\ell}} (z)=\dfrac{n_{\ell}c}{c-1}\mathcal{R}_{n_{\ell}}^{\beta_{\ell}}-\dfrac{c}{\left(c - 1\right)^2} n_{\ell} (\beta_{\ell} + n_{\ell}- 1)\mathcal{R}_{n_{\ell}-1}^{\beta_{\ell}}.
\end{equation}
To check equation \eqref{simple_eq_nablas_Rs} one can proceed as follows:
\begin{align*}
\nabla^{n_{\ell}}\rho^{\beta_{\ell}+n_{\ell}} (z)&=\nabla^{n_{\ell}-1}\left(\nabla\rho^{\beta_{\ell}+n_{\ell}} (z)\right)=\nabla^{n_{\ell}-1}\left(\rho^{\beta_{\ell}+n_{\ell}}(z)-\dfrac{\Gamma(z^r-1+\beta_{\ell}+n_{\ell})c^{z^{r}-1}}{\Gamma(z^r)}\right)\\
&=\nabla^{n_{\ell}-1}\left(\rho^{\beta_{\ell}+n_{\ell}}(z)-\dfrac{\Gamma(z^r-1+\beta_{\ell}+n_{\ell})c^{z^{r}}}{\Gamma(z^r+1)}\dfrac{(z^r+(\beta_{\ell}+n_{\ell}-1)-(\beta_{\ell}+n_{\ell}-1))}{c}
\right)\\
&=\left(\dfrac{c-1}{c}\right)\nabla^{n_{\ell}-1}\rho^{\beta_{\ell}+n_{\ell}}(z)+\dfrac{\beta_{\ell}+n_{\ell}-1}{c}\nabla^{n_{\ell}-1}\rho^{\beta_{\ell}+n_{\ell}-1}(z).\end{align*}
Here we used formula $\Gamma(z^r+1)=z^r\Gamma(z^r+1)$. Moreover, we added and subtracted the term $(\beta_{\ell}+n_{\ell}-1)$ to get the desired expression. Hence,
$$
\nabla^{n_{\ell}-1}\rho^{\beta_{\ell}+n_{\ell}}(z)=\left(\dfrac{c}{c-1}\right)\nabla^{n_{\ell}}\rho^{\beta_{\ell}+n_{\ell}} (z)-\dfrac{\beta_{\ell}+n_{\ell}-1}{c-1}\nabla^{n_{\ell}-1}\rho^{\beta_{\ell}+n_{\ell}-1}(z).
$$
Multiplying this equation from the left by $\left( \dfrac{c}{c-1} \right)^{n_{\ell}}\dfrac{n_{\ell}}{\rho^{\beta_{\ell}} (z)}$ gives \eqref{simple_eq_nablas_Rs}.

Finally, taking into account the following algebraic manipulation
\begin{align*}
\mathcal{R}_{n_{\ell}+1}^{\beta_{\ell}}&=\left( \dfrac{c}{c-1} \right)^{n_{\ell}+1}\dfrac{1}{\rho^{\beta_{\ell}} (z)}\nabla^{n_{\ell}}\left(\nabla\rho^{\beta_{\ell}+n_{\ell}+1} (z)\right)\\
&=\left( \dfrac{c}{c-1} \right)^{n_{\ell}+1}\dfrac{1}{\rho^{\beta_{\ell}} (z)}\nabla^{n_{\ell}}\left(\rho^{\beta_{\ell}+n_{\ell}+1}(z)-\dfrac{\Gamma(z^r+\beta_{\ell}+n_{\ell})c^{z^{r}-1}}{\Gamma(z^r)}\right)\\
&=\left( \dfrac{c}{c-1} \right)^{n_{\ell}+1}\dfrac{1}{\rho^{\beta_{\ell}} (z)}\nabla^{n_{\ell}}\left(\rho^{\beta_{\ell}+n_{\ell}} (z)\left((\beta_{\ell}+n_{\ell})+\dfrac{z^r(c-1)}{c}\right)\right)\\
&=\dfrac{c(\beta_{\ell}+n_{\ell})}{c-1}\mathcal{R}_{n_{\ell}}^{\beta_{\ell}}+\mathcal{R}_{n_{\ell}}^{\beta_{\ell}}z^r,
\end{align*}
and replacing the term $\mathcal{R}_{n_{\ell}}^{\beta_{\ell}}z^r$ with the expression \eqref{important_Rn} one obtains \eqref{eq_lemma_auxi}.

\end{proof}

\begin{remark}
Observe that formula \eqref{eq_lemma_auxi} gives a recurrence relation involving the polynomials with multi-indices $(n_{\ell}+1)\vec{e}_{\ell}$, $n_{\ell}\vec{e}_{\ell}$, and $(n_{\ell}-1)\vec{e}_{\ell}$ (see \eqref{eq:m1_rodrigues_proof_one_ray}), that is,
$$
M^{\vec{\beta},c}_{n_{\ell}\vec{e}_{\ell}+\vec{e}_{\ell}} (z)=z^{r}M^{\vec{\beta},c}_{n_{\ell}\vec{e}_{\ell}} (z)+ \left(\dfrac{c\left(n_{\ell}+\beta_{\ell}\right)}{c -1}+\dfrac{n_{\ell}}{c -1 }\right)M^{\vec{\beta},c}_{n_{\ell}\vec{e}_{\ell}} (z)
-\dfrac{c}{\left(c - 1\right)^2} n_{\ell} (\beta_{\ell} + n_{\ell}- 1)M^{\vec{\beta},c}_{n_{\ell}\vec{e}_{\ell}-\vec{e}_{\ell}} (z).
$$
\end{remark}

\begin{theorem} The type II Meixner-Angelesco polynomials of the second kind verify the following recurrence relations
\begin{equation}\label{eq:m1_recurrence_relation}
z^r M^{\vec{\beta},c}_{\vec{n}}(z)= M_{\vec{n}+\vec{e}_\ell}^{\vec{\beta},c}(z)+ b_{\vec{n},\ell}M^{\vec{\beta},c}_{\vec{n}}(z) + \sum_{j=0}^{r-1} d_{\vec{n}, j} M_{\vec{n}-\vec{e}_j}^{\vec{\beta},c}(z), \ \ \ \ell=0,\ldots,r-1,
\end{equation}
where $$b_{\vec{n},\ell} = \dfrac{c\left(n_{\ell}+\beta_{\ell}\right)+\lvert \vec{n}\rvert}{1 -c }\quad \text{and}\quad d_{\vec{n}, j} = \dfrac{c}{\left(c - 1\right)^2} n_j (\beta_{j} + n_{j}- 1).$$
\end{theorem}

\begin{proof} The action of the raising operator $\mathcal{R}_{n_{j}}^{\beta_{j}}$ on equation \eqref{eq_lemma_auxi} gives
\begin{equation}\label{eq_auxi_coro-1}
\mathcal{R}_{n_{j}}^{\beta_{j}} \mathcal{R}_{n_{\ell}+1}^{\beta_{\ell}}=\mathcal{R}_{n_{j}}^{\beta_{j}}\left(z^{r}\mathcal{R}_{n_{\ell}}^{\beta_{\ell}}\right)+ \left(\dfrac{c\left(n_{\ell}+\beta_{\ell}\right)}{c -1}+\dfrac{n_{\ell}}{c -1 }\right)\mathcal{R}_{n_{j}}^{\beta_{j}}\mathcal{R}_{n_{\ell}}^{\beta_{\ell}}
-\dfrac{c}{\left(c - 1\right)^2} n_{\ell} (\beta_{\ell} + n_{\ell}- 1)\mathcal{R}_{n_{j}}^{\beta_{j}}\mathcal{R}_{n_{\ell}-1}^{\beta_{\ell}},
\end{equation}
where $j,\ell\in\{0,\ldots,r-1\}$, $\ell\not=j$. Observe that from the multiplication-by-$z^r$ operator in \eqref{eq_lemma_auxi} and the commutative property of the operator \eqref{commutation_R_m-k} as well as its linearity one has
\begin{align}\label{eq_fundamental}
\mathcal{R}_{n_{j}}^{\beta_{j}}\left(z^{r}\mathcal{R}_{n_{\ell}}^{\beta_{\ell}}\right)
&=\left(\mathcal{R}_{n_{\ell}+1}^{\beta_{\ell}}- \dfrac{c\left(n_{\ell}+\beta_{\ell}\right)+n_{\ell}}{c -1}\mathcal{R}_{n_{\ell}}^{\beta_{\ell}}+\dfrac{c\,n_{\ell}}{\left(c - 1\right)^2}  (\beta_{\ell} + n_{\ell}- 1)\mathcal{R}_{n_{\ell}-1}^{\beta_{\ell}}\right)\mathcal{R}_{n_{j}}^{\beta_{j}}=z^r\mathcal{R}_{n_{\ell}}^{\beta_{\ell}}\mathcal{R}_{n_{j}}^{\beta_{j}}\notag\\
&=z^r\mathcal{R}_{n_{j}}^{\beta_{j}}\mathcal{R}_{n_{\ell}}^{\beta_{\ell}}\notag\\
&=\mathcal{R}_{n_{j}+1}^{\beta_{j}}\mathcal{R}_{n_{\ell}}^{\beta_{\ell}}-\left( \dfrac{c\left(n_{j}+\beta_{j}\right)+n_{j}}{c -1}\right)\mathcal{R}_{n_{j}}^{\beta_{j}}\mathcal{R}_{n_{\ell}}^{\beta_{\ell}}+\dfrac{c\,n_{j}}{\left(c - 1\right)^2}  (\beta_{j} + n_{j}- 1)\mathcal{R}_{n_{j}-1}^{\beta_{j}}\mathcal{R}_{n_{\ell}}^{\beta_{\ell}}\notag\\
&=z^r\mathcal{R}^{\beta_{j}}_{n_{j}}\mathcal{R}^{\beta_{\ell}}_{n_{\ell}}-\dfrac{n_{j}}{c -1 }\mathcal{R}_{n_{j}}^{\beta_{j}}\mathcal{R}_{n_{\ell}}^{\beta_{\ell}}+\dfrac{c}{\left(c - 1\right)^2} n_{j} (\beta_{j} + n_{j}- 1)\mathcal{R}_{n_{j}-1}^{\beta_{j}}\mathcal{R}_{n_{\ell}}^{\beta_{\ell}}.
\end{align}
Here, based on the property \eqref{operator_comm_L1} and expression \eqref{eq:m1_kth_times_operator_Psi}, the following algebraic manipulation was used
\begin{align*}
\mathcal{R}_{n_{j}+1}^{\beta_{j}}\mathcal{R}_{n_{\ell}}^{\beta_{\ell}}&=
\left(\dfrac{c}{c-1}\right)^{n_{j}+1} \left( \dfrac{1}{\rho^{\beta_{j}}(z)}\nabla^{n_{j}+1} \rho^{\beta_{j}+n_{j}+1}(z)\right)\mathcal{R}_{n_{\ell}}^{\beta_{\ell}}\\
&=\left(\dfrac{c}{c-1}\right)\left( \dfrac{1}{\rho^{\beta_{j}+n_{j}}(z)}\nabla \rho^{\beta_{\ell}+n_{j}+1}(z)\right)\mathcal{R}^{\beta_{j}}_{n_{j}}\mathcal{R}^{\beta_{\ell}}_{n_{\ell}}\\
&=\left(z^r+\dfrac{c(\beta_{j}+n_{j})}{c -1 }\right)\mathcal{R}^{\beta_{j}}_{n_{j}}\mathcal{R}^{\beta_{\ell}}_{n_{\ell}}.
\end{align*}
Therefore, \eqref{eq_auxi_coro-1} can be rewritten as follows
\begin{align}\label{recurrence_2_indices}
\mathcal{R}_{n_{j}}^{\beta_{j}} \mathcal{R}_{n_{\ell}+1}^{\beta_{\ell}}&=z^r\mathcal{R}_{n_{j}}^{\beta_{j}}\mathcal{R}_{n_{\ell}}^{\beta_{\ell}}
+ \left(\dfrac{c\left(n_{\ell}+\beta_{\ell}\right)}{c -1}+\dfrac{n_{\ell}+n_{j}}{c -1 }\right)\mathcal{R}_{n_{j}}^{\beta_{j}}\mathcal{R}_{n_{\ell}}^{\beta_{\ell}}-\dfrac{c}{\left(c - 1\right)^2} n_{j} (\beta_{j} + n_{j}- 1)\mathcal{R}_{n_{j-1}}^{\beta_{j}}\mathcal{R}_{n_{\ell}}^{\beta_{\ell}}\notag\\
&-\dfrac{c}{\left(c - 1\right)^2} n_{\ell} (\beta_{\ell} + n_{\ell}- 1)\mathcal{R}_{n_{j}}^{\beta_{j}}\mathcal{R}_{n_{\ell}-1}^{\beta_{\ell}}.
\end{align}
This expression gives the recurrence relation \eqref{eq:m1_recurrence_relation} involving the polynomial given in \eqref{eq:m1_rodrigues_proof_second_ray_equiv} (see also \eqref{polymomial-2-indices}) i.e.
\begin{align*}
z^r M^{\vec{\beta},c}_{n_{j}\vec{e}_{j}+n_{\ell}\vec{e}_{\ell}}(z)&= M_{n_{j}\vec{e}_{j}+n_{\ell}\vec{e}_{\ell}+\vec{e}_\ell}^{\vec{\beta},c}(z)+ 
\left(\dfrac{c\left(n_{\ell}+\beta_{\ell}\right)}{c -1}+\dfrac{n_{\ell}+n_{j}}{c -1 }\right)M^{\vec{\beta},c}_{n_{j}\vec{e}_{j}+n_{\ell}\vec{e}_{\ell}}(z)\\
& + \dfrac{c\,n_{\ell} (\beta_{\ell} + n_{\ell}- 1)}{\left(c - 1\right)^2} M_{n_{j}\vec{e}_{j}+n_{\ell}\vec{e}_{\ell}-\vec{e}_\ell}^{\vec{\beta},c}(z)
+ \dfrac{c\,n_{j} (\beta_{j} + n_{j}- 1)}{\left(c - 1\right)^2}M_{n_{j}\vec{e}_{j}-\vec{e}_j+n_{\ell}\vec{e}_{\ell}}^{\vec{\beta},c}(z).
\end{align*}

Finally, by iterating this process, i.e. by mathematical induction, in which the commutative linear operators \eqref{eq:m1_kth_times_operator_Psi}, $\mathcal{R}_{n_{k}}^{\beta_{k}}$, $k\in\{0,\ldots,r-1\}$, $k\not=j,\ell$ (see \eqref{commutation_R_m-k} as well as formula \eqref{eq_fundamental}) sequentially act on equation \eqref{recurrence_2_indices} one obtains relation \eqref{eq:m1_recurrence_relation}.
\end{proof}

\section{Conclusions}\label{conclusions}

In \cite{bib:arvesu-coussement-vanassche}, five families of discrete multiple orthogonal polynomials for AT-systems of measures were studied. Among these polynomials there are three families, namely, multiple Charlier and multiple Meixner polynomials of the first and second kind, respectively where the orthogonality measures are supported on an unbounded subset of $\mathbb{R}_{+}$. We refer to them as unbounded cases. The Angelesco systems with measures supported on unbounded rays of $\mathbb{C}$, with disjoint interior, which still lead to raising operators, Rodrigues formula, and recurrence relations among other algebraic properties is an interesting question. Indeed, it has not been studied until recently \cite{bib:arvesu-quintero1,bib:leurs-vanassche2,bib:leurs-vanassche1}. Before that few examples of Angelesco polynomials were known, namely, Jacobi-Angelesco, Jacobi-Laguerre, and Laguerre-Hermite polynomials. Only the last one considered two orthogonalizing weights with unbounded supports, that is, $\mathbb{R}_{-}$ and $\mathbb{R}_{+}$.

With the Meixner-Angelesco polynomials of the second kind we complete the study of Angelesco multiple orthogonal polynomials, which are analogues of the multiple Charlier and multiple Meixner polynomials of the first and second kind that have been studied for AT-system of discrete measures. We have obtained the raising operators, the Rodrigues-type formulas, the explicit expressions, and the nearest neighbor recurrence relations. These expressions involve complex parameters and are polynomials in $z^r$, which differ from \eqref{eq:m2_mclassic}.

Finally, other families of Angelesco polynomials should be studied. In particular, the bounded cases (Hahn-Angelesco and Kravchuk-Angelesco systems), where the orthogonality measures should be some complex analogues of the hypergeometric and binomial distributions supported on an $r$-star with bounded legs, respectively. In addition, some $q$-extensions involving Angelesco multiple orthogonal polynomials (in line with the AT-systems in \cite{bib:arvesu-qHahn,bib:arvesu-ramirez1,bib:arvesu-ramirez2,bib:arvesu-ramirez3}) should be investigated.

\section*{Acknowledgements} Jorge Arves\'u (J.A.) would like to thank the Department of Mathematics at Baylor University for hosting his visit in Spring 2021 which stimulated this research. The authors thank the anonymous referees for several helpful suggestions which improved this paper. The authors also thank Professor Andrei Mart\'{\i}nez-Finkelshtein, Baylor University, for his assistance and interesting discussions on the recurrence relations for Meixner-Angelesco multiple orthogonal polynomials.\\
\noindent{\bf Funding.} The research of J.A. was funded by Agencia Estatal de Investigaci\'on of Spain, grant number PID2021-122154NB-I00.

\end{document}